# Statistical approach based on 1D Voronoi tessellation as a tool for studying the randomness of fractional digits of some irrational numbers

M. Fanfoni
Dipartimento di Fisica Università di Roma Tor Vergata
Via della Ricerca Scientifica 00133 Roma Italy

M. Tomellini
Dipartimento di Scienze e Tecnologie Chimiche
Università di Roma Tor Vergata
Via della Ricerca Scientifica 00133 Roma Italy

## Abstract

An "experimental" study on the randomness of the fractional digits of $\pi$, $e$ and $\phi$ irrational numbers are presented. This is done by exploiting the *1D* Poisson-Voronoi tessellation. We employed two approaches and in both cases, within the numerical error, no differences have been detected between the irrational fractional digits and an equivalent random sequence of digits. The number of tested digits is $1.6 \times 10^7$ and $4 \times 10^7$ for the first and second approach, respectively. Although not shown here, we investigated several irrational numbers and all of them have displayed a similar behavior.

## 1-Introduction

The nucleation process, which is at the basis of new phase formation, gives rise to a geometrical pattern of nuclei which can be characterized by the Voronoi tessellation. Each nucleus, which can be approximated as a dimensionless point, is at the center of its own cell. On the other hand the spatial distribution of nuclei is usually correlated because of the diffusion process that rules the growth. Particularly interesting is the model case of a random distribution of nuclei, leading to the well known Poisson-Voronoi tessellation [1-3].

The tassellation can be characterized by the distribution function of cells' size (also called in nucleation theory "capture zone") $f(x)$ [4].

It is worth reminding that the exact form of $f(x)$ is known only for the *1D* case. Kiang, in fact, demonstrated that the distribution of the Voronoi's tiles is described by the single parameter Gamma distribution [5-8]

$$\gamma(x;\alpha)=\frac{\alpha^{\alpha}}{\Gamma(\alpha)}x^{\alpha-1}e^{-\alpha x} \qquad (1)$$

with $\alpha=2$ i.e. ,

$$\gamma(x;2)=4x\,e^{-2x}. \qquad (2)$$

The Voronoi tessellations have been often trotted out in relation to the spatial distribution of nuclei when they are formed on a substrate during film growth. In that case one has to do with a *2D* problem and, what is more, due to the monomer diffusion and/or the strain field, the interaction between nuclei switches on and, depending on its intensity, a certain degree of spatial correlation takes place making eqn.1 starts failing. Specifically, if the distribution of nuclei is not random the value of the $\alpha$ parameter increases with respect to that of the random distribution. It is important to point out that the $\alpha$ value is very sensitive to the degree of correlation among nuclei [9].

As a matter of fact, in case of Poissonian distribution of points in a *2D* space, the Voronoi tessellation provides a distribution which can be very well fitted by eqn.1 with $\alpha \in [3.5,3.6]$. As anticipated, to this day, an analytical demonstration of this result is still not available. On the other hand, random nucleation which takes place on a linear defect (step) leads to a *1D* tessellation which can be modeled through eqn.1

The randomness of the digits of irrational numbers is a long standing question, see for example [10-13] and references therein. Here we would like to give a contribution to the "debate" with a sort of "experimental" analysis of some irrational numbers: π, *e* (Nepero's number), ϕ (Golden ratio), γ (Eulero-Mascheroni constant), $e^{\pi}$ (Gelfond's constant), $\sqrt{2}$ , $\sqrt{3}$ , $ln2$ ,$2^{\sqrt{2}}$ (Gelfond-Schneider constant).
In order to do that, we exploit the exact solution for Poisson-Voronoi tessellation in *1D*. Specifically, given a *1D* lattice we conveniently (see below) select from the fractional digits of the irrational number under examination (from now on $\boldsymbol{I_n}$ ) the coordinates of the nuclei. Successively, we determine the *1D* Voronoi tessellation, fit eqn.1 to it having α as parameter and evaluate how much it differs from the exact value for the Poissonian distribution α=2.

## 2-Methodology

Two routines have been developed. As far as the first routine is concerned, let us explain in formal way what the code carries out. Let $S = \{1,2 \ldots M-1\}$ be a discrete set of integers. We define the map

$$F\colon \underbrace{S \times S \ldots \times S}_{h} \rightarrow \underbrace{S \times S \ldots \times S}_{h}$$

which works as follow: $h$ random integer numbers $(n_1, \ldots n_h)$ are chosen in such a way that $n_k \in S \; \forall \; k \in [1, h]$; $l(n_k)$ is the length of $n_k$, i.e. its number of digits (e.g. $l(3458)=4$).
If $\boldsymbol{D_d^{(0)}} = (a_1, a_2, a_3, \ldots)$ is the set of fractional digits of the irrational number $\boldsymbol{I_n} = A.a_1a_2a_3 \ldots$, $h$ nuclei are then selected whose coordinates are $(c_1(l(n_1)), \ldots, c_h(l(n_h)))$. By way of an example, let us consider the case of $\pi$ with M=10000 and $h$=6; $\boldsymbol{I_n} = 3.14159265358979323846\ldots$ and $\boldsymbol{D_d^{(0)}} = (1,4,1,5,9,2,6,5,3,5,8,9,7,9,3,2,3,8,4,6)$. To the random set $(8921, 6325, 1124, 199, 650, 81) \equiv (n_1, \ldots n_6)$ is associated the set of their lenght $((4,4,4,3,3,2) \equiv (l(n_1), \ldots l(n_6))$ . Accordingly, the set $(c_1, \ldots c_6)$ is defined as follows:

$$c_1 = a_1a_2a_3a_4 = 1415$$
$$c_2 = a_5a_6a_7a_8 = 9265$$
$$c_3 = a_9a_{10}a_{11}a_{12} = 3589$$
$$c_4 = a_{13}a_{14}a_{15} = 793$$
$$c_5 = a_{16}a_{17}a_{18} = 238$$
$$c_{46} = a_{19}a_{20} = 46$$

so, in the set of M-1 0's (i.e. *1D* lattice) the positions $(c_1, \ldots c_6)$ are transformed into 1's (nuclei).
As far as the coordinates of $c_k$ are concerned, the general case takes the form

$c_1 = a_1 a_2 \dots a_{l(n_1)}$
$c_2 = a_{l(n_1)+1} a_{l(n_1)+2} \cdots a_{l(n_1)+l(n_2)}$
$c_3 = a_{l(n_1)+l(n_2)+1} a_{l(n_1)+l(n_2)+2} \cdots a_{l(n_1)+l(n_2)+l(n_3)}$

and in compact notation it reads

$$c_k = a_{\sum_{i=1}^{k-1} l(n_i)+1} a_{\sum_{i=1}^{k-1} l(n_i)+2} \cdots a_{\sum_{i=1}^{k} l(n_i)}$$

The map, in short notation, is $F(n_1, \dots, n_k) = (c_1(n_1), \dots, c_p(n_h))$. In other words, the first $l(n_1)$ digits of $\boldsymbol{D}_{\boldsymbol{d}}^{(\mathbf{0})}$ give rise to the coordinate of the first nucleus on the *1D* lattice, these digits are ruled out from $\boldsymbol{D}_{\boldsymbol{d}}^{(\mathbf{0})}$and of the new set $\boldsymbol{D}_{\boldsymbol{d}}^{(\mathbf{1})} = \boldsymbol{D}_{\boldsymbol{d}}^{(\mathbf{0})} \backslash (a_1 a_2 \dots a_{l(n_1)})$ the first $l(n_2)$ digits give rise to the coordinate of the second nucleus. For example, if $l(n_1) = 4$, the coordinate of the first nucleus is $c_1 = a_1 a_2 a_3 a_4$ and $\boldsymbol{D}_{\boldsymbol{d}}^{(\mathbf{1})} = \boldsymbol{D}_{\boldsymbol{d}}^{(\mathbf{0})} \backslash (a_1, a_2, a_3, a_4) = (a_5, a_6, a_7, \dots)$.
The first routine determines $N_0$ coordinates in such a way to reach a given density $\rho = \frac{N_0}{M}$, and then determines the Voronoi cells distribution. Eqn.1is fitted to the size cell distribution and the parameter $\alpha$ determined. The second routine finds the set $\{p_k^i\}$ of the position of the $i-th$ digit $d_i$ in $\boldsymbol{D}_{\boldsymbol{d}}^{(\mathbf{0})}$. If $P_M^i = Max[p_k^i]$ is the maximum of the set $\{p_k^i\}$, one can consider a *1D* lattice of dimension $P_M^i$, where the positions $\{p_k^i\}$ plays the role of the coordinates of the nuclei which give rise to a Voronoi tesselation of the lattice. For example, the first 41 digits of $\pi$ are:

$\{1, 4, 1, 5, \mathbf{9}, 2, 6, 5, 3, 5, 8, \mathbf{9}, 7, \mathbf{9}, 3, 2, 3, 8, 4, 6, 2, 6, 4, 3, 3, 8, 3, 2, 7,$
$\mathbf{9}, 5, 0, 2, 8, 8, 4, 1, \mathbf{9}, 7, 1, 6\}.$

If one is interested, for instance, in the digit 9, the previous set transforms in
$\{0, 0, 0, 0, \mathbf{1}, 0, 0, 0, 0, 0, 0, \mathbf{1}, 0, \mathbf{1}, 0, 0, 0, 0, 0, 0, 0, 0, 0, 0, 0, 0, 0, 0, 0, \mathbf{1},$
$0, 0, 0, 0, 0, 0, 0, \mathbf{1}, 0, 0, 0\},$

where the 1's occupy the positions of the digit 9. Nevertheless, the lattice to be tessellated is taken to be

$$\{\mathbf{1}, 0, 0, 0, 0, 0, 0, \mathbf{1}, 0, \mathbf{1}, 0, 0, 0, 0, 0, 0, 0, 0, 0, 0, 0, 0, 0, 0, 0, \mathbf{1}, 0, 0, 0, 0, 0, 0, 0, \mathbf{1}\},$$

on the other hand the first and last 0's are removed.

Thus, once the distribution of the Voronoi segments' length is determined, it is fitted to eqn.1.
In this way, not only one can check the frequency of the ten $d_i$ digits, but also their "spatial" distribution within $\boldsymbol{D}_{\boldsymbol{d}}^{(\mathbf{0})}$.

## 3-Results

*3.1 First routine*

Several values of $Card(\boldsymbol{D}_{\boldsymbol{d}}^{(\mathbf{0})})$ have been employed obtaining, in fact, the same results. The data here presented are related to of $Card\left(\boldsymbol{D}_{\boldsymbol{d}}^{(\mathbf{0})}\right) = 1.6 \times 10^7$. The value of the density is $\rho = 10^{-3}$, while the lattice length is $M = 10^8$. $Card(\boldsymbol{D}_{\boldsymbol{d}}^{(\mathbf{0})})$ is much more inferior to the values which one finds in up to date articles (above all concerning $\pi$), our results show quite clearly that data stabilize pretty soon. It goes without saying that the fruition of a larger computer would allow to push the check well beyond the frontiers here reached.
Table 1 summarizes the final data for $\left|2 - \alpha_{I_n}\right|$ where $\alpha$ is the parameter of the Gamma distribution as obtained by non-linear fit routine of Mathematica9 by Wolfram Research.

**Table 1**

| Rnd | π | e | ϕ |
|---|---|---|---|
| 0.0022 | 0.0027 | 0.0022 | 0.0027 |

The standard error of the fits ranges from a minimum of 0.003 to a maximum of 0.005. In fig.1 the experimental distributions and the respective best fit of some irrational numbers ($\pi, e,$ ϕ) of Tab.1 are reported. Similar outputs are obtained for the other numbers.
In the hypothesis the fractional part of the irrational numbers were random (as it seems to be from the values of α's), it would be significant to calculate the mean $\alpha$ value: in our case it is ⟨α⟩=1.99841, which differs from α =2 of 0.0795%.
As a last test for the sensitivity of our first routine to randomness of the $\boldsymbol{D}_{\boldsymbol{d}}^{(\mathbf{0})}$ set we wrote another code which works as follows: in the first place, among the 10! possible permutations of the simple set $\{1,2,3,4,5,6,7,8,9,0\}$ , the code extract at random $10^3$ different permutations. The new object is a set of $10^3$ distinct elements each being a set of ten numbers. Among the $10^3!$ possible permutations, the code extracts 1600 permutations. As a result, a set of $10 \times 1000 \times 1600 = 1.6 \times 10^7$ numbers is obtained. The last set is then processed in the same way as the $1.6 \times 10^7$ digits of $\boldsymbol{I_n}$. In fig.2 the calculated distribution function is shown. The change is dramatic, it is only a remembrance of a Gamma distribution, it is quite close to a $\delta_+(x)$ function.

*3.2 Second routine*

As far as the second routine is concerned, the analysis has been restricted to π, *e*, ϕ and the random set. In this case, $Card\left(\boldsymbol{D_d^{(0)}}\right) = 4 \times 10^7$ and the average recurrence of any digits is, on average, $f = 4 \times 10^6$. In fig.3 the distribution of one digit is shown while Tab.2 gives the value of $|2 - \alpha\,|$'s the *1D* Voronoi distributions obtained from the positions of any single digit within the lattice made up of $4 \times 10^6$ for a random distribution, and the three irrational numbers π, *e*, and ϕ. Although we report here the results for a single digit, the analysis has been performed for all digits obtaining, in fact, similar outputs.

Let $S'$ be a set of *N* random numbers each of them chosen in the interval (0, 9). We are interested in determining the positions or coordinates of a given digit inside *S'* and subsequently its *1D* Voronoi tessellation and relative tile distribution function. It goes without saying that the probability that at the position k there is a given digit d is *p*=0.1, while the probability of its absence, $\bar{d}$ is

$q = 1 - p$. Accordingly, the string $d\bar{d}\bar{d}\bar{d}\bar{d}\bar{d}d\bar{d}\bar{d}\bar{d}\bar{d}\bar{d}\bar{d}d = d\bar{d}^{n_1}d\bar{d}^{n_2}d$ has probability to occur given by $p^2q^{n_1+n_2}$ . Let us define the new variable $s = \frac{n_1+n_2}{2} + 1$ which is the tile of the sequence, in other words the probability for the tile $s$ to occur is given by

$$P(s) = \sum_{k=1}^{2s-1} p^2q^{2(s-1)} = 2(s-1)p^2q^{2(s-1)}$$

The *P(s)* distribution is normalized, in fact, if $k = 2s - 1$ then $\sum_{k=1}^{\infty} P(s) = 1$. The *s* mean value is

$\langle s \rangle = \langle \frac{k+1}{2} \rangle = \sum_{k=1}^{\infty} \frac{k+1}{2} P(k) = \frac{1}{p}$, which is the expected result. In order to bridge the gap between the discrete and the continuum case we generate a hystogram $\{s - 1/2, P(s)\}$ where $s \in [\frac{1}{2}, 50]$ with step 0.05. The hystogram has been normalized and properly rescaled in such a way that$\langle s \rangle = 1$. Then the Gamma distribution eqn.1 has been fitted to it receiving α=2.0021. On this basis, if the distribution of the digits of an irrational number, $\boldsymbol{I_n}$, is random, the recurrence of the position of each digit within the $\boldsymbol{D_d^{(0)}}$ set is expected to be in agreement with $\gamma(x; 2)$

**Table 2**

| Rnd | $\pi^{[d=3]}$ | $e^{[d=2]}$ | $\phi^{[d=1]}$ |
|---|---|---|---|
| 0.0021 | 0.0011 | 0.0032 | 0.0013 |

In Tab.2 the digit ($n$) indicated in the first row is the one whose distribution has been displayed in fig.3. Comparable results have been obtained for the other digits. Noteworthy also in this case, as we done for the data of Tab.1, the difference between the average of the $\alpha_{I_n}$, and α=2, is found to be of the order of magnitude of 0.09%.

In order to test the sensitivity of the present approach to distinguish between a random and a non-random set of nuclei, we performed numerical simulation for a correlated set of nuclei. Specifically, the distribution of nuclei is chosen in accord to the hard core model where the distance among nuclei cannot be shorter than a given value, say $l^*$. The quantity $\rho l^*$ is therefore a measure of the displacement of the system from the random case ; it is related to the fraction of length forbidden to nucleation.

We perform computer simulations of space tessellation for nucleation process with increasing $l^*$ , determine the *P(s)* distribution and, eventually, the α parameter through fit of eqn.1 to the *P(s)*. The results of this computation are reported in fig.4 which show a quite linear trend of α with $l^*$. In fact, the curve is in excellent accord with the second order polinomial $(\alpha - 2) = al^* + bl^{*2}$ with coefficients α=0.0047 and $b = 1.042 \times 10^{-3}$ . This shows the high sensitivity of α with the randomness of the system; the lowest correlation $l^* = 1$ would imply an α variation larger than the deviation attained in Tab.2.

**4- Conclusions**

In conclusion, the two routines here developed seem to establish that the fractional digits of irrational numbers (at least those here analyzed) are random distributed. The present approach exploits the description of a physical process which implies random events as previously presented in the literature in connection to the diffusion (random walk).

**Acknowledgement**

The authors are indebted with Prof. S. Trapani for the helpful discussions.

## Refereneces

## Figures

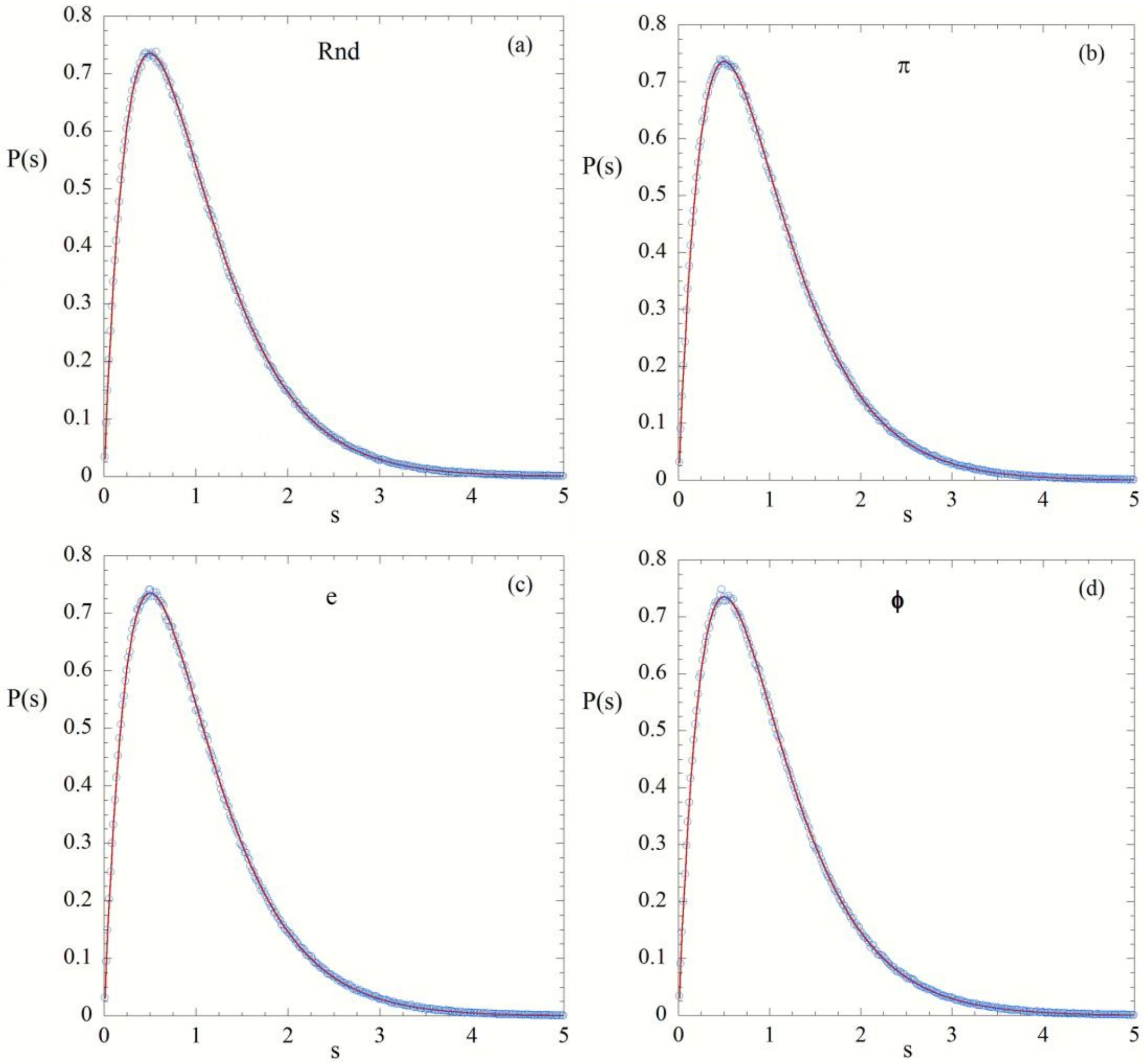

**Fig.1** Cell-size Voronoi distribution as obtained from the first routine as described in the text. Here, we report the outputs for the random sequence (panel a), irrational numbers $\pi$ (panel b), e (panel c) and ϕ (panel d) as solid symbols. The solid line is the best fit of eqn.1 to the data points. The values of fitting parameter $\left|2-\alpha_{I_n}\right|$ are displayed in Tab.1, together with those of other irrational numbers.

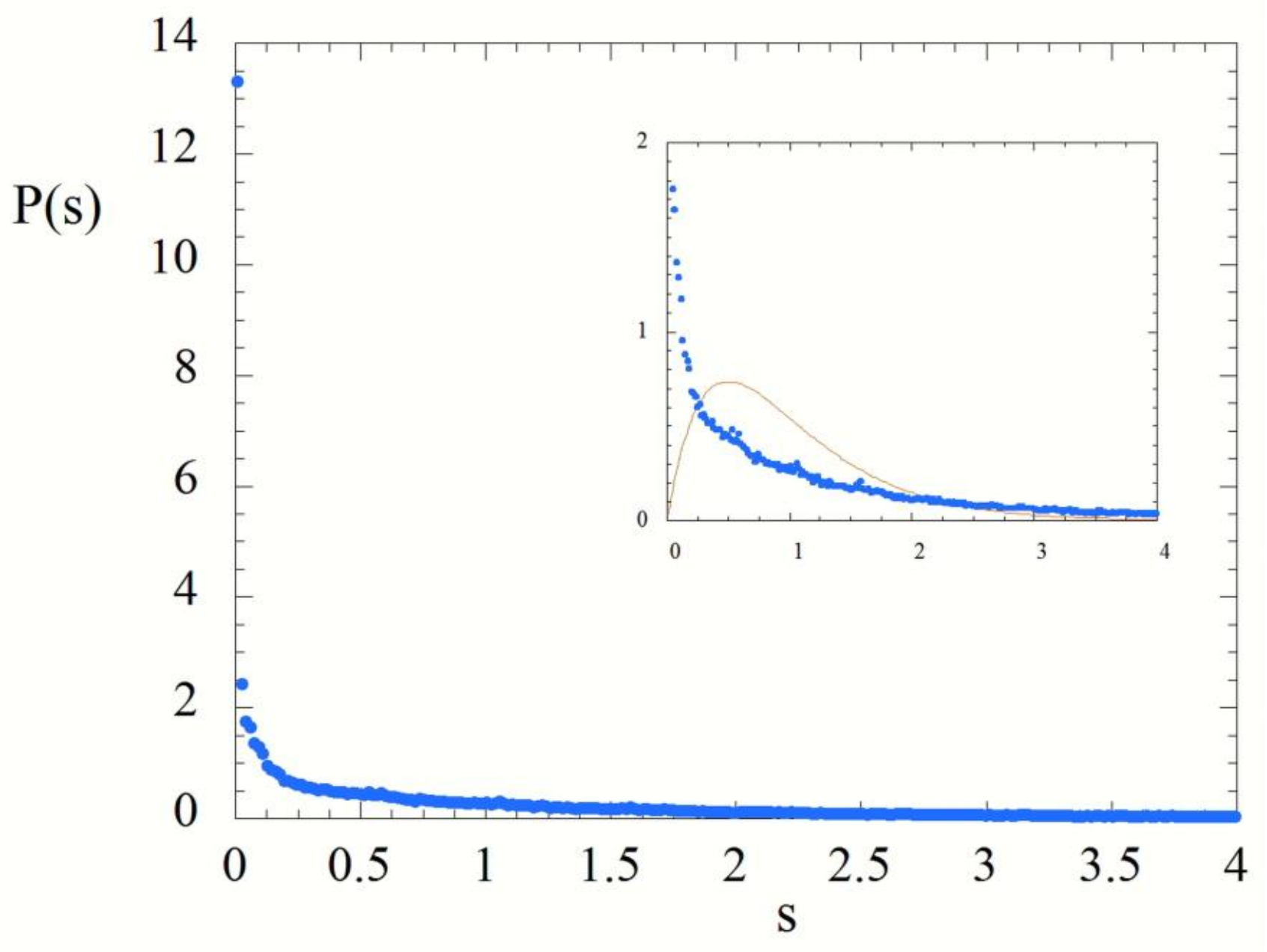

**Fig.2** Voronoi distribution for the "test" on randomness as described in the text. The first routine is applied to a sequence of $1.6 \times 10^7$ digits chosen according to the procedure reported in the text. In the inset it is evidenced the comparison with the $\gamma(x; 2)$distribution (solid line)

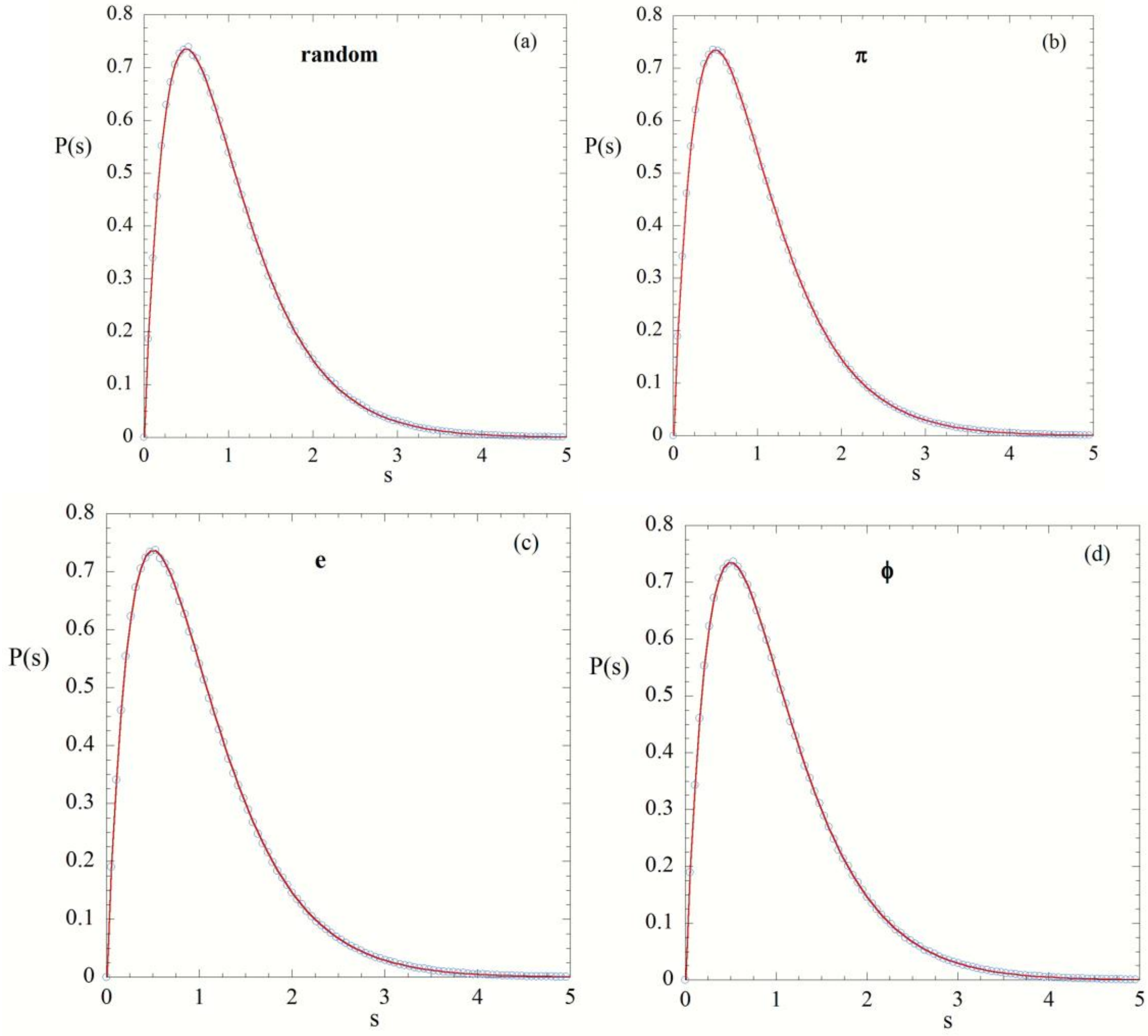

**Fig.3** Cell-size Voronoi distribution as obtained from the second rutine as described in the text. Here, we report the outputs for the random sequence (panel a), irrational numbers $\pi$, (digit 3 panel b), e (digit 2 panel c), and ϕ (digit 1 panel d) as solid symbols. The solid line is the best fit of eqn.1 to the data points. The values of fitting parameter $\left|2-\alpha_{I_n}\right|$ are displayed in Tab.2, where $d$ denotes the considered digit.

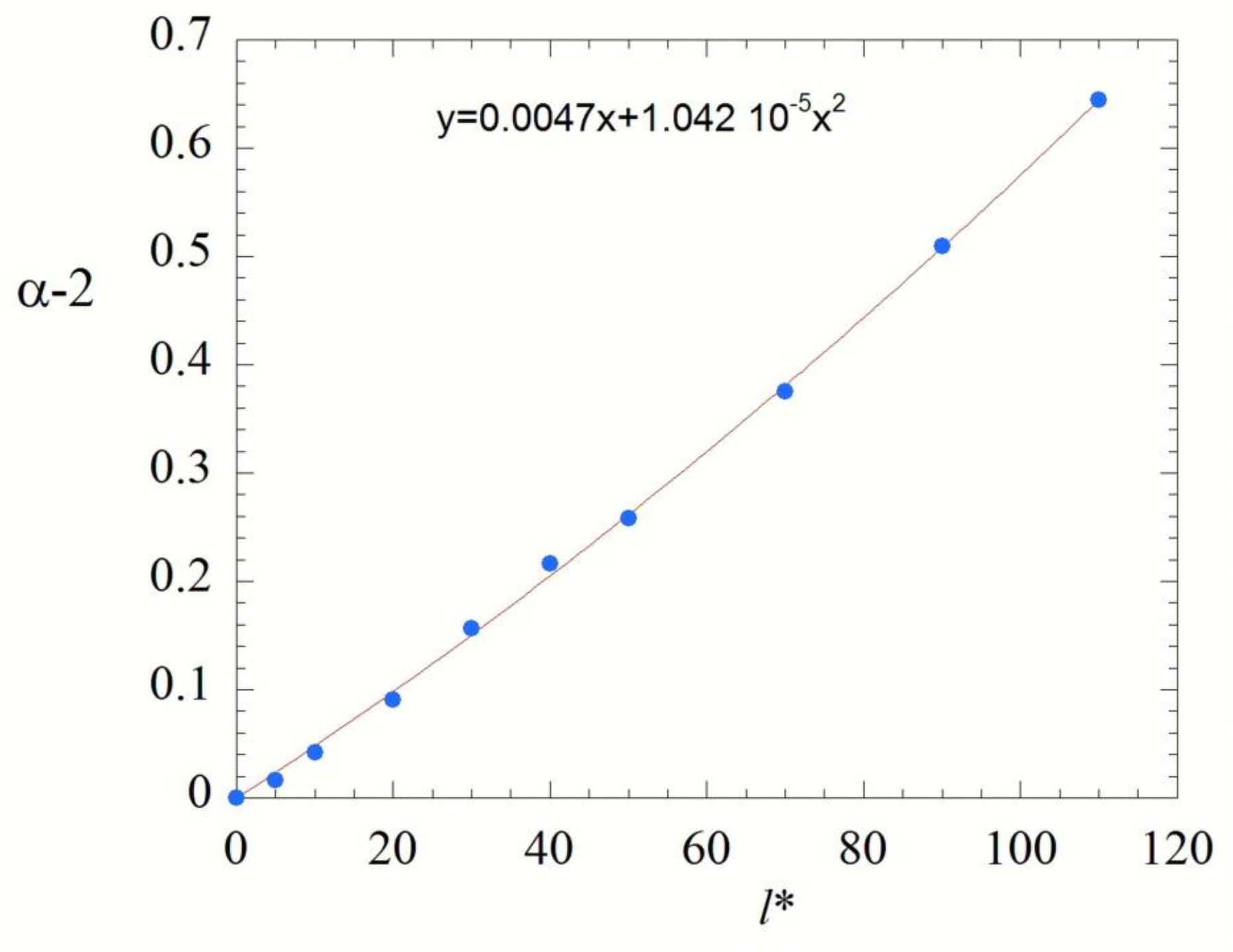

**Fig.4** Behavior of the a parameter as a function of the correlation degree among nuclei, $l^*$. Solid line is the best fit of the second order polynomial to the data points. Fit parameters evidence the high sensitivity of a to correlation degree.